\documentclass[11pt,a4paper]{article}
\usepackage[T1]{fontenc}
\usepackage{lmodern}
\usepackage[margin=28mm]{geometry}
\usepackage{amsmath,amssymb,amsthm,microtype}
\usepackage[colorlinks=true,linkcolor=blue]{hyperref}
\usepackage{orcidlink}

\newtheorem{maintheorem}{Theorem}

\newtheorem{lemma}{Lemma}[section]

\newcommand{\tes}{\operatorname{tes}}
\title{Total Edge Irregularity Strength of Cubic Graphs}
\author{Aleams Barra$^1$\,\orcidlink{0000-0001-8618-4821}, \qquad Suhadi Wido Saputro$^2$\,\orcidlink{0000-0003-3886-400X}, \qquad Hajime Tanaka$^3$\,\orcidlink{0000-0002-5958-0375} \\[3mm]
\small $^1$Algebra Research Group, Faculty of Mathematics and Natural Sciences, \\ \small Institut Teknologi Bandung, Indonesia \\ \small\texttt{aleamsbarra@itb.ac.id} \\
\small $^2$Combinatorial Mathematics Research Group, Faculty of Mathematics and Natural Sciences, \\ \small Institut Teknologi Bandung, Indonesia \\ \small\texttt{suhadi@itb.ac.id} \\
\small $^3$Research Center for Pure and Applied Mathematics \\
\small Graduate School of Information Sciences, Tohoku University, Japan \\
\small \texttt{htanaka@tohoku.ac.jp}}
\date{}
\begin{document}
\maketitle
\begin{abstract}
We show that every finite simple cubic graph on $n$ vertices has total edge irregularity strength $n/2+1$, confirming the Ivan\v{c}o--Jendrol' conjecture for cubic graphs.
The key step is to algorithmically construct a $\{1,n/2+1\}$-vertex labeling such that half of the vertices are labeled $1$ and the average degree of the subgraph induced by them is close enough to two.
We extend such a labeling to a total edge irregular $(n/2+1)$-labeling, possibly after a local adjustment to a single vertex.
\end{abstract}

\section{Introduction}

Let $s$ be a positive integer.
A \emph{total edge irregular $s$-labeling} of a finite simple graph $G$ with vertex set $V(G)$ and edge set $E(G)$ is a mapping
$\varphi:V(G)\cup E(G)\to[1,s]:=\{1,2,\ldots,s\}$ such that the edge weights
$\varphi(u)+\varphi(uv)+\varphi(v)$ are distinct for all $uv\in E(G)$.
The \emph{total edge irregularity strength} of $G$, denoted by $\tes(G)$,
is the least positive integer $s$ for which such a labeling exists. This concept was introduced by Ba\v{c}a, Jendrol', Miller, and Ryan~\cite{BJMR2007DM} as a total-labeling analogue of the classical irregularity strength of a graph, introduced by Chartrand and collaborators~\cite{chartrand}.

For every graph $G$ with a nonempty edge set, Ba\v{c}a et al.~\cite{BJMR2007DM} established the bounds
\begin{equation}\label{lower bound}
    \tes(G) \geq \max\left\{ \left\lceil\frac{{|E(G)|}+2}{3}\right\rceil,\left\lceil \frac{\Delta(G)+1}{2} \right\rceil\right\},
\end{equation}
where $\Delta(G)$ is the maximum degree of $G$. They also determined the exact values for several elementary graph families, including paths, cycles, wheels, and friendship graphs. Subsequently, Ivan\v{c}o and Jendrol'~\cite{ivancojendrol}
conjectured that this lower bound is attained by every graph
except $K_5$. More precisely, they conjectured that
\[
\tes(G)=\max\left\{
\left\lceil\frac{|E(G)|+2}{3}\right\rceil,
\left\lceil\frac{\Delta(G)+1}{2}\right\rceil
\right\},
\]
which they subsequently proved to hold for all trees.

Since its introduction, determining $\operatorname{tes}(G)$ has remained an active topic of investigation. Exact values have been established for a wide variety of graph families, including complete graphs and complete bipartite graphs~\cite{complete}, fan and helm graphs and other wheel-related graphs~\cite{jeyanthi}, and various graph products~\cite{strongpaths,categoricalpaths,cycleproduct,prism}.
The subject has remained active in recent years, with further results
for generalized arithmetic staircase graphs~\cite{stair},
Cartesian products of bipartite graphs and paths~\cite{wijaya},
several classes of cycle-related graphs~\cite{cycle2024},
and arithmetic ladder graphs~\cite{ladder}. These examples are by no means exhaustive; further results and references can be found in Gallian's periodically updated survey of graph labeling~\cite{gallian}. The continuing expansion of this literature reflects the sustained interest in total edge irregularity strength and related labeling problems.

Probabilistic methods have also played an important role in obtaining exact values of $\tes(G)$ for general graphs. Brandt, Mi\v{s}kuf, and Rautenbach~\cite{bmr08} proved that $\tes(G)=\lceil(|E(G)|+2)/3\rceil$ under a sufficient condition relating the maximum degree, order, and size of $G$. They subsequently strengthened their result~\cite{bmr09}, showing that the same equality holds whenever $|E(G)|>111000\Delta(G)$. A major further advance was made by Pfender~\cite{pfender}, who established the Ivan\v{c}o--Jendrol' conjecture for every graph with at least $7\times10^{10}$ edges. Thus, probabilistic techniques have established the conjectured formula for broad classes of graphs and, ultimately, for all sufficiently large graphs.

Despite Pfender's general result, considerable attention has continued to be devoted to determining the total edge irregularity strength of particular families of cubic graphs. Exact results are known for generalized Petersen graphs~\cite{haque}, toroidal fullerenes~\cite{fullerene}, disjoint unions of prisms~\cite{disjointprisms}, and disjoint unions of $(3,6)$- and $(4,6)$-fullerene graphs~\cite{fullerene36,fullerene46}. Such results remain valuable even in light of Pfender's sweeping theorem because they often provide explicit, deterministic constructions of optimal labelings, rather than solely probabilistic existence arguments. Moreover, these constructions typically apply without a lower bound on the order of the graphs, thereby complementing the general large-graph theorem and demonstrating how the structures of particular cubic graph families can be exploited to attain the conjectured lower bound.

In the present paper, we completely resolve the problem of determining the total edge irregularity strength of cubic graphs. We show that every finite simple cubic graph attains the conjectured lower bound, independently of its order and whether it is connected. Our result therefore encompasses both connected and disconnected cubic graphs and extends the known exact evaluations for particular families to the entire class. The main theorem is as follows.

\begin{maintheorem}\label{thm:main}
Every finite simple cubic graph $G$ on $n$ vertices satisfies
\[
\tes(G)=\frac n2+1.
\]
\end{maintheorem}

We organize the proof of Theorem~\ref{thm:main} as follows. In Section~\ref{sec:admissibility}, we establish the \emph{admissibility lemma}
(Lemma~\ref{lem:admissible}) for general finite simple graphs, which
characterizes when a vertex labeling using only the labels $1$ and $s$
can be extended to a total edge irregular $s$-labeling.
We then prove the main theorem in Section~\ref{sec:proof}.
The key idea is to construct a vertex labeling
in which nearly half of the vertices receive the label $1$, and the subgraph
induced by these vertices has average degree close to $2$.
We obtain such a $\{1,s\}$-vertex labeling algorithmically.
In one exceptional case, the resulting labeling does not satisfy
the conditions of the admissibility lemma, but this can be resolved by adjusting
the label of a single vertex.

\section{Admissibility lemma}\label{sec:admissibility}

The admissibility criterion underlying our approach was established by Ivan\v{c}o and Jendrol'~\cite[Lemma~3]{ivancojendrol}, who used binary vertex labelings to construct total edge irregular labelings. We independently arrived at the same sufficient conditions through the parallel admissibility criterion of Barra and Afifurrahman~\cite[Lemma~2.1]{barraafifurrahman}, which concerns completing a $\{1,s\}$-edge labeling to a vertex irregular total labeling. We subsequently became aware of Ivan\v{c}o and Jendrol's earlier result.
In the following, we present their criterion in our notation as an if-and-only-if characterization and give a direct constructive proof for completeness.

For a vertex labeling $\varphi$ of a simple graph $G$, put
\[
E_r=\{uv\in E(G):\varphi(u)+\varphi(v)=r\} \qquad (r=2,3,\dots).
\]

\begin{lemma}[Admissibility]\label{lem:admissible}
Let $s\ge3$, and let $\varphi:V(G)\to\{1,s\}$ be a vertex labeling of a finite simple graph $G$.
Write
$(a,b,c)=(|E_2|,|E_{s+1}|,|E_{2s}|)$.
Then, $\varphi$ extends to a total edge irregular $s$-labeling of $G$ if and only if
\begin{equation}\label{admissibility}
a,b,c\le s,\qquad a+b,b+c\le2s-1, \qquad a+b+c\,(={|E(G)|})\le3s-2.
\end{equation}
\end{lemma}
\begin{proof}
For necessity, note that the weights of the edges in $E_2$, $E_{s+1}$, and $E_{2s}$ lie in $[3,s+2]$, $[s+2,2s+1]$, and $[2s+1,3s]$, respectively. Each interval has size $s$, the union of two consecutive intervals has size $2s-1$, and the union of all three has size $3s-2$. Distinct edge weights therefore require \eqref{admissibility}.

For sufficiency, we assign consecutive edge labels within each class, choosing the smallest starting label that places its weight interval strictly to the right of the preceding one. Put
\[
\alpha=\max\{1,a+2-s\},\qquad
\beta=\max\{1,\alpha+b+1-s\}.
\]
Assign the edges in $E_2$, $E_{s+1}$, and $E_{2s}$ distinct labels within their respective classes from
\[
[1,a],\qquad[\alpha,\alpha+b-1],\qquad[\beta,\beta+c-1],
\]
respectively. All edge labels are positive. Their largest labels satisfy $a\le s$ and
\[
\begin{aligned}
\alpha+b-1&=\max\{b,a+b+1-s\}\le s,\\
\beta+c-1&=\max\{c,b+c+1-s,a+b+c+2-2s\}\le s,
\end{aligned}
\]
by \eqref{admissibility}.
Thus, every edge label lies in $[1,s]$.

The resulting edge weights are distinct within each class and form the intervals
\begin{equation*}
    [3,a+2], \qquad [s+1+\alpha,s+\alpha+b], \qquad [2s+\beta,2s+\beta+c-1].
\end{equation*}
Moreover, since
\begin{gather*}
    s+1+\alpha\ge s+1+(a+2-s)=a+3>a+2, \\
    2s+\beta\ge 2s+(\alpha+b+1-s)=s+\alpha+b+1>s+\alpha+b,
\end{gather*}
these three intervals are ordered from left to right and are pairwise disjoint.
Hence, all edge weights are distinct, completing the proof.
\end{proof}

\section{Proof of Theorem~\ref{thm:main}}\label{sec:proof}

Let $G$ be a finite simple cubic graph on $n$ vertices, and put $s=n/2+1\,(\ge 3)$.
In view of \eqref{lower bound}, it suffices to find a total edge irregular $s$-labeling of $G$.

Note that the vertex labelings $\varphi:V(G)\to\{1,s\}$ are in bijection with the subsets $U$ of $V(G)$ by letting $U=\varphi^{-1}(1)$.
We retain the notation in the previous section, and recall $(a,b,c)=(|E_2|,|E_{s+1}|,|E_{2s}|)$, where $a+b+c={|E(G)|}=3n/2=3s-3$.
Since $G$ is cubic, it follows that $3|U|=2a+b$ and $3|V(G)\setminus U|=2c+b$.

In the following, we will work with subsets $U\subset V(G)$ with $|U|=n/2$.
Note that $a=c$ in this case.
We will begin with a subset $U$ such that $a\ge s-2$.
Whenever $a>s$, our algorithm will allow us to reduce $a$ by three by replacing a vertex in $U$.

We choose the initial subset $U\subset V(G)$ as follows.
Let $G_1,G_2,\ldots,G_\ell$ be the connected components of $G$, and write $n_i=|V(G_i)|$.
Let $k\ge0$ denote the largest integer such that
\[
    \rho:=\sum_{i=1}^k n_i<n/2,
\]
where an empty sum is zero.
Note that $n_{k+1}\ge n/2-\rho$.
Let $T$ be a spanning tree of $G_{k+1}$, and pick any subtree $S$ of $T$ having $n/2-\rho$ vertices.
Let
\[
    U=V(G_1)\cup\cdots\cup V(G_k)\cup V(S).
\]
We have $|U|=n/2$.
Moreover, since each $G_i$ is cubic and $S$ is a tree, we have
\begin{equation*}
    a\ge \sum_{i=1}^k|E(G_i)|+|E(S)|=\frac{3\rho}{2}+\left(\frac{n}{2}-\rho-1\right)\geq \frac{n}{2}-1=s-2,
\end{equation*}
as desired.

Assume now that $a\,(=c)>s$.
The average degrees of the subgraphs of $G$ induced by $U$ and $V(G)\setminus U$ are equal and are given by $\overline{d}_U$, where
\begin{equation*}
    \overline{d}_U=\frac{2a}{n/2}=\frac{2a}{s-1}>\frac{2s}{s-1}>2.
\end{equation*}
It follows that we can find a vertex $u\in U$ (resp.~$v\in V(G)\setminus U$) all three of whose neighbors are in $U$ (resp.~$V(G)\setminus U$).
Note that $u$ and $v$ are not adjacent.
We update $U$ as follows:
\[
    U\leftarrow(U\setminus\{u\})\cup\{v\}.
\]
Observe that $|U|$ remains $n/2$, and that the integer $a$ is then updated as
\begin{equation*}
    a \leftarrow a-3.
\end{equation*}

Repeat the above procedure until we first get $a\le s$.
We then have $a\in\{s-2,s-1,s\}$, so that the final distribution $(a,b,c)$ is one of
\[
    (s,s-3,s),\qquad(s-1,s-1,s-1),\qquad(s-2,s+1,s-2).
\]
The first two cases above satisfy \eqref{admissibility} and thus extend to total edge irregular $s$-labelings by Lemma~\ref{lem:admissible}.
Hence, for the rest of the proof, we assume that $(a,b,c)=(s-2,s+1,s-2)$.

Since $\overline{d}_U=2(s-2)/(s-1)<2$ in this case, we can find a vertex $u\in U$ that has at most one neighbor in $U$.
Change the label of $u$ from $1$ to $2$ (so the corresponding vertex labeling $\varphi$ now takes three values $1,2$, and $s$).
If $u$ has no neighbor in $U$, then we have
\begin{equation*}
    (|E_2|,|E_{s+1}|,|E_{s+2}|,|E_{2s}|)
=(s-2,s-2,3,s-2),
\end{equation*}
and we complete the edge labels as follows:
\[
\begin{array}{c|c|c}
\hline\hline
\text{Class}&\text{Edge labels}&\text{Edge weights}\\ \hline
E_2 &[1,s-2]&[3,s]\\
E_{s+1}&[1,s-2]&[s+2,2s-1]\\
E_{s+2}&[s-2,s]&[2s,2s+2]\\
E_{2s}&[3,s]&[2s+3,3s] \\
\hline\hline
\end{array}
\]
On the other hand, if $u$ has exactly one neighbor in $U$, then we have
\begin{equation*}
    (|E_2|,|E_3|,|E_{s+1}|,|E_{s+2}|,|E_{2s}|)
=(s-3,1,s-1,2,s-2),
\end{equation*}
and we complete the edge labels as follows:
\[
\begin{array}{c|c|c}
\hline\hline
\text{Class}&\text{Edge labels}&\text{Edge weights}\\ \hline
E_2 &[1,s-3]&[3,s-1]\\
E_3&\{s-2\}&\{s+1\}\\
E_{s+1}&[1,s-1]&[s+2,2s]\\
E_{s+2}&[s-1,s]&[2s+1,2s+2]\\
E_{2s}&[3,s]&[2s+3,3s] \\
\hline\hline
\end{array}
\]
(We note that $E_2=\emptyset$ if $s=3$, i.e., $G=K_4$.)
In each case, all labels are in $[1,s]$, and the edge weight sets are pairwise disjoint.
Hence, we obtain a total edge irregular $s$-labeling.
This completes the proof of Theorem \ref{thm:main}.

\section*{Acknowledgments}
AB and SWS are supported by the ITB Research Program 2026 under the ITB International Research Scheme through the Directorate of Research and Innovation, Institut Teknologi Bandung (Project ID: FMIPA.PN-6-124-2026). HT was supported by JSPS KAKENHI Grant Number \href{https://kaken.nii.ac.jp/en/grant/KAKENHI-PROJECT-26K06717/}{JP26K06717}.

The authors acknowledge the use of ChatGPT for language editing and polishing of the Introduction. All suggestions were reviewed and revised by the authors, who retain full responsibility for the entire text.

\end{document}